\documentclass[11pt,reqno]{amsart}
\usepackage[T1]{fontenc}
\usepackage[utf8]{inputenc}
\usepackage{lmodern}
\usepackage{amsmath,amssymb,amsthm,mathtools,mathrsfs}
\usepackage[a4paper,margin=28mm]{geometry}
\usepackage[expansion=false]{microtype}
\usepackage{booktabs,array}
\usepackage{hyperref}

\numberwithin{equation}{section}
\newtheorem{theorem}{Theorem}[section]
\newtheorem{proposition}[theorem]{Proposition}
\newtheorem{lemma}[theorem]{Lemma}
\newtheorem{corollary}[theorem]{Corollary}
\theoremstyle{remark}
\newtheorem{remark}[theorem]{Remark}
\title[Chow groups and odd homology of Chow varieties]{Chow groups and odd homology of Chow varieties of space curves}
\author[Y. Chen]{Youming Chen }
\address{School of Mathematical Science, Chongqing University of Technology, Chongqing 400054, P.R. China}%
\email{youmingchen@cqut.edu.cn}%

\author[W. Hu]{Wenchuan Hu}
\address{School of Mathematics, Sichuan University, Chengdu 610064, P.R. China}%
\email{huwenchuan@gmail.com}%

\date{\today}


\begin{document}

\begin{abstract}
Let $X_d=C_{1,d}(\mathbb{P}^{3})$ denote the Chow variety of effective curves of degree $d$ in complex projective space $\mathbb{P}^{3}$.
For every $d\geq13$, we prove that the rational vector space of algebraically trivial classes
in $\operatorname{CH}_{4d-37}(X_d)$ is uncountably dimensional.
We also establish that $H^{8d-37}(X_d,\mathbb{Q})$ has nonzero
Hodge numbers with respect to the weight filtration in bidegrees $(4d,4d-37)$ and $(4d-37,4d)$,
so that $b_{8d-37}(X_d)\geq2$.
The construction uses a component of the locus of nonspecial genus-$10$ space curves,
together with a canonical form on a smooth projective model of the universal Picard variety.
An extremal Hodge-type argument extends the resulting correspondence to the entire Chow variety,
which may be non-equidimensional.
These results disprove our Chow group comparison conjectures 
and the  conjectured vanishing of ordinary odd Betti numbers.
We also show that finite affine pavings fail to exist
and that the specified degree-stabilization maps are not injective.
\end{abstract}

\maketitle

\section{Introduction}

Throughout this paper, all varieties are defined over $\mathbb{C}$.
A variety is a reduced separated scheme of finite type; in particular, it need not be irreducible or equidimensional.
Homology and cohomology are those of the associated analytic space, with rational coefficients unless otherwise stated.
For a projective variety $X$, we write
$$
b_j(X)=\dim_{\mathbb{Q}}H_j(X,\mathbb{Q})
      =\dim_{\mathbb{Q}}H^j(X,\mathbb{Q}).
$$
These are ordinary singular Betti numbers.
They also equal the ranks of the corresponding integral homology groups,
since a complex projective variety has the homotopy type of a finite CW complex.

For integers $n\geq1$, $0\leq p\leq n$, and $d\geq1$,
let $C_{p,d}(\mathbb{P}^{n})$ denote the Chow variety parametrizing effective algebraic $p$-cycles of degree $d$ in complex projective $n$-space.
These varieties are projective, but in general they are singular, reducible, and not equidimensional.
Their topology in a stable range is known, whereas their Chow groups and their homology outside this range are less well understood.  These questions and  several comparison conjectures are surveyed in  \cite{Hu2021}.

Let $\operatorname{CH}_{q}(X)$ denote the group of algebraic cycles of dimension $q$ modulo rational equivalence.
Unadorned Chow groups have integral coefficients.
For an abelian group $A$,
set $A_{\mathbb{Q}}=A\otimes_{\mathbb{Z}}\mathbb{Q}$, and write
$$
\operatorname{CH}_{q}(X)_{\mathrm{alg}}
 =\{\alpha\in\operatorname{CH}_{q}(X):\alpha\sim_{\mathrm{alg}}0\}.
$$
Rational and algebraic equivalence coincide on $q$-cycles precisely when this subgroup  vanishes.
We use Chow groups and their functorial operations in the sense of Fulton~\cite{Fulton}.
The rational cycle class map is denoted by
$$
\operatorname{cl}_{q,X,\mathbb{Q}}:
\operatorname{CH}_{q}(X)_{\mathbb{Q}}\longrightarrow H_{2q}(X,\mathbb{Q}).
$$
In particular, none of the arguments below requires an intersection ring on a singular target.

The second author proposed that the Chow groups of $C_{p,d}(\mathbb{P}^{n})$
are isomorphic to  the corresponding even degree integral singular homology groups \cite[Conj.~4.9]{Hu2021}.
We restated this assertion and its rational version in~\cite[Conj.~1.1 and~1.3]{ChenHuChow}.
 They also proposed the vanishing of all ordinary odd Betti numbers;
 see the paragraph following Conjecture~1.3 in either of~\cite{ChenHuChow,ChenHuLawson}.
 These assertions go beyond the established stable range.
 The rational cycle class map is an isomorphism for $0\leq q\leq d$ by~\cite[Cor.~2.4]{ChenHuChow},
 while the homotopy stability theorem gives
\begin{equation*}
H_j\bigl(C_{1,d}(\mathbb{P}^{3}),\mathbb{Q}\bigr)=0
\quad\text{for odd $j$ with $j\leq2d$};
\end{equation*}
see~\cite[Thm.~1 and Cor.~5]{Hu2015}.

We study the Chow variety of space curves and set
\begin{equation*}
X_d=C_{1,d}(\mathbb{P}^{3})\quad(d\geq1),\quad
q_d=4d-37,\quad k_d=8d-37\quad(d\geq13).
\end{equation*}
Our first result concerns cycles \emph{on} $X_d$, rather than the curves in $\mathbb{P}^{3}$ parameterized by $X_d$.
As usual, ``very general'' means outside a countable union of proper Zariski closed subsets of an irreducible parameter space.

\begin{theorem}\label{thm:chow}
For every integer $d\geq13$, the rational vector space
$$
\operatorname{CH}_{q_d}(X_d)_{\mathrm{alg}}
 \otimes_{\mathbb{Z}}\mathbb{Q}
$$
has uncountable dimension.
More precisely, there exists an irreducible variety $B^{\circ}$
and an algebraic family of effective $q_d$-cycles $F_b$ on $X_d$, indexed by $b\in B^{\circ}(\mathbb{C})$,
such that all the $F_b$ are algebraically equivalent and their rational Chow classes take uncountably many distinct values.
For every fixed $b_0\in B^{\circ}(\mathbb{C})$, a very general $b$ satisfies
$$
[F_b]-[F_{b_0}]\neq0
\quad\text{in }\operatorname{CH}_{q_d}(X_d)_{\mathbb{Q}}.
$$
In particular, these differences give non-torsion integral classes that are algebraically trivial but not rationally trivial.
\end{theorem}

Our second result concerns ordinary odd homology and its mixed Hodge structure.
For a rational mixed Hodge structure $H$, let $W_{\bullet}$ and $F^{\bullet}$ denote its weight and Hodge filtrations,
and set
$$
h^{a,b}(H)=\dim_{\mathbb{C}}
 \operatorname{Gr}^{a}_{F}\operatorname{Gr}^{W}_{a+b}
 (H\otimes_{\mathbb{Q}}\mathbb{C}).
$$
Here $\operatorname{Gr}^{W}_{w}=W_w/W_{w-1}$ and $\operatorname{Gr}^{a}_{F}=F^a/F^{a+1}$.
We use Deligne's mixed Hodge structures~\cite{DeligneII,DeligneIII} on algebraic varieties.

\begin{theorem}\label{thm:odd}
For every integer $d\geq13$, one has
\begin{equation}\label{eq:main-hodge}
\begin{split}
h^{4d,\,4d-37}\!\left(\operatorname{Gr}^{W}_{k_d}H^{k_d}(X_d,\mathbb{Q})\right)&\geq1,\\
h^{4d-37,\,4d}\!\left(\operatorname{Gr}^{W}_{k_d}H^{k_d}(X_d,\mathbb{Q})\right)&\geq1.
\end{split}
\end{equation}
Consequently,
$$
b_{k_d}(X_d)=\operatorname{rank}_{\mathbb{Z}}H_{k_d}(X_d,\mathbb{Z})\geq2.
$$
\end{theorem}

Since $q_d>d$ and $k_d>2d$ for $d\geq13$, these results lie outside the established stable ranges.
The algebraically trivial differences in Theorem~\ref{thm:chow} have zero class in $H_{2q_d}(X_d,\mathbb{Q})$.
The odd homology in Theorem~\ref{thm:odd} instead detects the correspondence defined by the whole family.
In particular, this does not assert that an algebraically trivial cycle has a nonzero ordinary cycle class.

We record the numerical example and the main consequences separately.
Their proofs are given in Section~\ref{sec:applications}.

\begin{corollary}\label{cor:24}
For $(d,q)=(13,15)$ or $(24,59)$, the vector space
$$
\operatorname{CH}_{q}\bigl(C_{1,d}(\mathbb{P}^{3})\bigr)_{\mathrm{alg}}
 \otimes_{\mathbb{Z}}\mathbb{Q}
$$
has uncountable dimension.
In particular, rational and algebraic equivalence do not coincide for $59$-dimensional cycles on $C_{1,24}(\mathbb{P}^{3})$.
Moreover,
$$
b_{67}(X_{13})\geq2,\quad
b_{155}(X_{24})\geq2.
$$
\end{corollary}

\begin{corollary}\label{cor:comparisons}
For every $d\geq13$, the kernel of
$$
\operatorname{cl}_{q_d,X_d,\mathbb{Q}}:
\operatorname{CH}_{q_d}(X_d)_{\mathbb{Q}}
 \longrightarrow H_{2q_d}(X_d,\mathbb{Q})
$$
has uncountable dimension.
The two vector spaces in this display are not isomorphic even abstractly.
 Hence the Chow group assertions of~\cite[Conj.~4.9]{Hu2021} and~\cite[Conj.~1.1 and~1.3]{ChenHuChow} are false.
\end{corollary}

\begin{corollary}\label{cor:paving}
For every $d\geq13$, the mixed Hodge structure $H^{k_d}(X_d,\mathbb{Q})$ is not Hodge--Tate,
and $X_d$ admits no finite affine paving.
\end{corollary}

A mixed Hodge structure is called Hodge--Tate
if the Hodge numbers of its graded pieces of the weight filtration vanish outside the diagonal bidegrees $(a,a)$.
A finite affine paving is a finite filtration by closed subsets whose successive differences are affine spaces.
These notions concern ordinary cohomology, not the virtual Hodge polynomial.

For the stabilization consequence, fix a line $\ell\subset\mathbb{P}^{3}$ and, for integers $e\geq d\geq1$,
define the proper cycle addition morphism
\begin{equation*}
s_{d,e}:X_d\longrightarrow X_e,\quad
Z\longmapsto Z+(e-d)[\ell].
\end{equation*}
The next corollary uses the rational comparison theorem~\cite[Cor.~2.4]{ChenHuChow}.
This result is not used in the proofs of Theorems~\ref{thm:chow} and~\ref{thm:odd}.

\begin{corollary}\label{cor:stabilization}
Let $d\geq13$ and $e\geq q_d$. Then
$$
(s_{d,e})_*\bigl(\operatorname{CH}_{q_d}(X_d)_{\mathrm{alg}}
 \otimes_{\mathbb{Z}}\mathbb{Q}\bigr)=0.
$$
Consequently the kernel of
$$
(s_{d,e})_*:\operatorname{CH}_{q_d}(X_d)_{\mathbb{Q}}
 \longrightarrow\operatorname{CH}_{q_d}(X_e)_{\mathbb{Q}}
$$
has uncountable dimension. In particular, addition of $35[\ell]$ gives a map
$$
\operatorname{CH}_{59}(X_{24})_{\mathbb{Q}}
 \longrightarrow\operatorname{CH}_{59}(X_{59})_{\mathbb{Q}}
$$
with an uncountable-dimensional kernel.
\end{corollary}

The common Hodge-theoretic argument also has a consequence for arbitrary projective varieties.
For a smooth projective variety $Y$, let $\Omega_Y^r$ denote the sheaf of regular differential $r$-forms and set
$$
h^{r,0}(Y)=\dim_{\mathbb{C}}H^0(Y,\Omega_Y^r).
$$
The following general statement is a consequence of Proposition~\ref{prop:extreme};
it does not depend on the construction of the particular varieties $X_d$.

\begin{corollary}\label{cor:general-hodge}
Let $V$ be an irreducible component of a projective variety $X$, let $\dim V=N$,
and let $Y\longrightarrow V$ be a projective resolution of singularities. For every $1\leq r\leq N$, set $k=2N-r$.
Then
\begin{equation}\label{eq:extreme-bound}
h^{N,N-r}\!\left(\operatorname{Gr}^{W}_{k}H^k(X,\mathbb{Q})\right)
 \geq h^{r,0}(Y),\quad
b_k(X)\geq2h^{r,0}(Y).
\end{equation}
No equidimensionality assumption on $X$ is required.
\end{corollary}

Finally, for a projective variety $X$ and $j\geq2r\geq0$, let
$$
L_rH_j(X)=\pi_{j-2r}\bigl(\mathcal{Z}_r(X)\bigr),
$$
where $\mathcal{Z}_r(X)$ is the topological group of algebraic $r$-cycles with its Lawson topology,
and $\pi_i$ denotes the $i$-th homotopy group based at the zero-cycle.
Theorem~\ref{thm:odd} and the Lawson comparison in~\cite[Prop.~4.8]{Hu2021} give the following consequence.

\begin{corollary}\label{cor:lawson}
For every $d\geq13$, one has
\begin{equation*}
\dim_{\mathbb{Q}}L_1H_{k_d}(X_d)_{\mathbb{Q}}\geq2.
\end{equation*}
\end{corollary}

Theorem~\ref{thm:odd} disproves the conjectured vanishing of ordinary odd Betti numbers in~\cite{ChenHuChow,ChenHuLawson}.
It does not by itself refute the virtual-Hodge conjecture~\cite[Conj.~4.1]{Hu2021}.
Corollary~\ref{cor:lawson} is an application of a cited Lawson comparison, not a counterexample to it.
Nor do the results address the symmetric product conjecture~\cite[Conj.~1.2]{ChenHuChow}
or assert that the constructed classes remain nonzero modulo algebraic equivalence.

The proofs of the two main theorems share the same geometric construction.
A smooth projective model of the universal Picard variety of genus-$10$ curves has dimension $37$.
The divisor formulas of Logan and Farkas--Popa, together with the integral Picard group of the stack of stable pointed curves
and the extension theorem of Farkas--Verra, produce a nonzero canonical form in Picard degree $10$.
We check the symmetric group character and use a trace between Picard torsors to transfer this form to every Picard degree.
In particular, no restriction $3\nmid d$ is imposed.
Nonspecial very ample line bundles then give a genuine $4d$-dimensional component $V_d\subset X_d$ dominating the degree-$d$ Picard space,
with general geometric fibre of dimension $q_d$.

An extremal Hodge type survives pushforward from a resolution of $V_d$ to the whole Chow variety.
It gives Theorem~\ref{thm:odd} directly. Combined with the Bloch--Srinivas method,
the same surviving class prevents the fibre correspondence from being rationally constant, yielding Theorem~\ref{thm:chow}. Section~\ref{sec:correspondences} establishes the general arguments and proves Corollary~\ref{cor:general-hodge}.
Sections~\ref{sec:picard} and~\ref{sec:components} construct the canonical forms and the Chow component, respectively. Section~\ref{sec:applications} proves the two main theorems and the remaining corollaries.
We make no claim about the minimal degree or cycle dimension in which these phenomena occur.

\section{Components and cycle correspondences}\label{sec:correspondences}

\subsection{Extremal Hodge types}

We retain the Hodge-number and filtration conventions from the introduction.
For a smooth variety $Y$, write $\omega_Y=\Omega_Y^{\dim Y}$ for its canonical line bundle.
For a coherent sheaf $\mathcal{F}$ on a projective variety $Y$, set $h^i(Y,\mathcal{F})=\dim_{\mathbb{C}}H^i(Y,\mathcal{F})$.
When $Y$ is smooth and projective, $H^{a,b}(Y)=H^b(Y,\Omega_Y^a)$
denotes the corresponding Hodge summand of $H^{a+b}(Y,\mathbb{C})$.

The Hodge indices of the cohomology of a proper variety of dimension at most $s$ lie between $0$ and $s$.
Dually, the indices of its homology lie between $-s$ and $0$.
This dimension bound, the strictness of morphisms,
and the exactness of the mixed Hodge structures on the localization sequences are standard consequences of~\cite[Sec.~8]{DeligneIII}.
We apply them to varieties that need not be equidimensional.

\begin{proposition}\label{prop:extreme}
Let $V$ be an irreducible component of a projective variety $X$, let $\dim V=N$,
and let
$$
j:Y\longrightarrow V\hookrightarrow X
$$
be a projective resolution of $V$. For $1\leq r\leq N$, set $k=2N-r$.
Then
\begin{equation}\label{eq:extreme-injection}
H^0(Y,\Omega_Y^r)\longrightarrow H_k(X,\mathbb{C}),
\quad \theta\longmapsto j_*(\theta\cap[Y]),
\end{equation}
is injective.
\end{proposition}

\begin{proof}
Choose a nonempty open subset $U\subset V$ over which $j$ is an isomorphism,
and remove its intersections with the other irreducible components of $X$.
The resulting $U$ is also open in $X$. Set $E=Y\setminus j^{-1}(U)$.
Then $\dim E\leq N-1$, and the localization long exact sequence for Borel--Moore homologies gives an exact sequence of mixed Hodge structures
\begin{equation}\label{eq:BM}
H_k(E,\mathbb{Q})\longrightarrow H_k(Y,\mathbb{Q})
\longrightarrow H_k^{\mathrm{BM}}(U,\mathbb{Q}).
\end{equation}
Here $H_k^{\mathrm{BM}}$ denotes Borel--Moore homology.
Both $E$ and $Y$ are proper, so their Borel--Moore and ordinary homology agree.

For $\theta\in H^{r,0}(Y)$, Poincar\'e duality identifies $\theta\cap[Y]$
with a class of Hodge type $(r-N,-N)$ in the pure Hodge structure $H_k(Y)$.
The image of $H_k(E)$ has no such Hodge type, since both of its indices are at least $-(N-1)$.
Strictness in~\eqref{eq:BM} therefore shows that the restriction of such a nonzero class to $U$ is nonzero.
Proper pushforward and restriction commute, and $j$ is an isomorphism over $U$.
This proves~\eqref{eq:extreme-injection}.
\end{proof}

\begin{proof}[Proof of Corollary~\ref{cor:general-hodge}]
Set $k=2N-r$ and apply Proposition~\ref{prop:extreme}.
The image of $H_k(Y)$ under $j_*$ is a pure Hodge structure of weight $-k$.
Its Hodge subspace of type $(r-N,-N)$ has dimension at least $h^{r,0}(Y)$.
Duality between $H_k(X,\mathbb{Q})$ and $H^k(X,\mathbb{Q})$,
followed by Hodge symmetry, gives the first inequality in~\eqref{eq:extreme-bound}.
Hodge symmetry on $\operatorname{Gr}^{W}_{k}H^k(X,\mathbb{Q})$ gives the same lower bound for the complementary type $(N-r,N)$.
Since $r>0$, their dimensions add, proving the Betti number inequality.
\end{proof}

The component assumption is essential.
For an arbitrary proper subvariety $V\subset X$,
a dense open subset of $V$ need not be open in $X$, and the last step of the proof would not apply.
Proposition~\ref{prop:extreme} makes no assertion about injectivity on all of $H_k(Y)$.

\subsection{Rational equivalence in a family}

Let $B$ be a smooth projective variety of dimension $m$ and let $X$ be projective.
A class $\Gamma\in\operatorname{CH}_{m+q}(B\times X)$ defines
$$
\Gamma_*:\operatorname{CH}_0(B)\longrightarrow\operatorname{CH}_q(X)
$$
by refined intersection in the smooth factor $B$ followed by proper pushforward.
Its action on singular cohomology and homology is
\begin{equation}\label{eq:correspondence-action}
\Gamma_{\#}(\omega)
=(p_X)_*\bigl(p_B^*(\omega)\cap\operatorname{cl}(\Gamma)\bigr),
\quad
\Gamma_{\#}:H^m(B)\longrightarrow H_{m+2q}(X),
\end{equation}
where $p_B,p_X$ are the projections and $\operatorname{cl}(\Gamma)$ is the singular homology class of $\Gamma$.
Both constructions remain valid  when $X$ is singular.

The following form of the Bloch--Srinivas argument~\cite{BlochSrinivas} is necessary here.
We include the passage to the generic fibre since smoothness of the target is not assumed.

\begin{lemma}\label{lem:spread}
Suppose that $B$ is connected, and let $B^{\circ}\subset B$
be a nonempty open subset over which a representative of $\Gamma\in\operatorname{CH}_{m+q}(B\times X)$
forms a flat family of $q$-cycles.
If the set
$$
\{[\Gamma_b]:b\in B^{\circ}(\mathbb{C})\}
\subset\operatorname{CH}_q(X)_{\mathbb{Q}}
$$
is countable, then
\begin{equation}\label{eq:spread}
\Gamma=B\times Z+\Gamma_D
\quad\text{in }\operatorname{CH}_{m+q}(B\times X)_{\mathbb{Q}},
\end{equation}
where $Z\in\operatorname{CH}_q(X)_{\mathbb{Q}}$, $D\subsetneq B$ is closed,
and $\Gamma_D$ is supported on $D\times X$.
\end{lemma}

\begin{proof}
Choose points $b_i\in B^{\circ}(\mathbb{C})$ representing the countably many rational Chow classes.
Rational equivalences between effective cycles are parameterized by countably many finite type spaces:
add an auxiliary effective cycle to both endpoints and use rational chains in a Chow space, fixing their degrees.
One may use two-pointed genus-zero stable maps to this projective Chow space,
whose properness does not require a smooth target~\cite{FultonPandharipande}.
Signed cycles are treated by moving negative terms to the other side;
positive multiples account for rational coefficients.
This is the usual parameter space argument for rational equivalence;
compare~\cite[Prop.~8 and Cor.~9]{BanerjeeGuletskii} and~\cite[Chap.~I, Sec.~3]{Kollar}.

For each $i$ and each positive integer $a$,
these parameter spaces describe the condition $a(\Gamma_b-\Gamma_{b_i})\sim_{\mathrm{rat}}0$.
Their constructible images cover $B^{\circ}(\mathbb{C})$.
An irreducible variety over $\mathbb{C}$ is not a countable union of proper closed subvarieties.
Consequently, one of the parameter spaces dominates $B^{\circ}$.
Let $K=\mathbb{C}(B)$
and choose a closed point in its nonempty generic fibre.
After a finite field extension $K'/K$, and a further finite extension if needed for a universal cycle datum,
we obtain
$$
\Gamma_{K'}=(\Gamma_{b_i})_{K'}
\quad\text{in }\operatorname{CH}_q(X_{K'})_{\mathbb{Q}}.
$$
Restriction followed by proper pushforward is multiplication by $[K':K]$,
so this equality already holds over $K$.
A rational equivalence consists of finitely many cycles and rational functions,
and thus extends  over a nonempty open subset of $B$.
Chow localization for $B\times X$ gives~\eqref{eq:spread},
with $Z=\Gamma_{b_i}$ and with $D$ enlarged to include the excluded closed subset.
The argument uses cycle spaces and localization, not an intersection product on $X$.
\end{proof}

\subsection{Fibre correspondences}

\begin{theorem}\label{thm:fibres}
Let $V$ be an $N$-dimensional irreducible component of a projective variety $X$.
Suppose that there exists a dominant rational map $V\dashrightarrow B$,
where $B$ is a smooth connected projective variety of dimension $m>0$ and $H^0(B,\Omega_B^m)\neq0$.
Set $q=N-m$. Then
$\dim_{\mathbb{Q}}\bigl(\operatorname{CH}_q(X)_{\mathrm{alg}}
\otimes_{\mathbb{Z}}\mathbb{Q}\bigr)$
is uncountable.
A resolution of the graph defines a family of effective $q$-cycles over an irreducible base whose fibre differences realize this conclusion.
For every fixed parameter $b_0$ in a suitable flat open subset,
a very general $b$ gives a difference which is nonzero over $\mathbb{Q}$.
\end{theorem}

\begin{proof}
Resolve the graph to obtain a smooth projective $Y$ and morphisms
$$
f:Y\longrightarrow B,\quad j:Y\longrightarrow V\hookrightarrow X,
$$
where $j$ is birational to $V$ and $f$ is dominant. Set
$$
\Gamma=(f,j)_*[Y]\in\operatorname{CH}_N(B\times X).
$$
Pullback of regular differential forms by a dominant rational map between smooth projective varieties in characteristic zero is injective.
We use the usual birational invariance of $H^0(\Omega^r)$,
so the choice of resolution does not change this space.
Proposition~\ref{prop:extreme} and the projection formula show that
\begin{equation}\label{eq:nonzero-action}
\Gamma_{\#}(\omega)=j_*\bigl(f^*\omega\cap[Y]\bigr).
\end{equation}
This action is injective on $H^{m,0}(B)$.

After shrinking to a nonempty open $B^{\circ}$, the fibres of $f$ give effective $q$-cycles $F_b$ on $X$ by proper pushforward.
These cycles are algebraically equivalent.
Indeed, the zero-cycles of degree zero on the connected smooth projective $B$ are algebraically trivial,
and the correspondence $\Gamma$ preserves algebraic equivalence.

Suppose that the classes of the $F_b$ in $\operatorname{CH}_q(X)_{\mathbb{Q}}$ formed a countable set.
Lemma~\ref{lem:spread} would give~\eqref{eq:spread}.
Both terms on the right-hand side act trivially on $H^{m,0}(B)$.
For $B\times Z$, integration over the first factor is zero in this degree.
For $\Gamma_D$, resolve each subvariety in its support.
Its projection to $B$ has image of dimension less than $m$, so every holomorphic $m$-form pulls back to zero. Equation~\eqref{eq:correspondence-action} then gives zero.
This contradicts~\eqref{eq:nonzero-action}.

Thus the differences $[F_b]-[F_{b_0}]$ assume uncountably many values in the algebraically trivial subgroup tensored with $\mathbb{Q}$.
A vector space of countable Hamel dimension over the countable field $\mathbb{Q}$ is itself countable.
Hence their span has uncountable dimension.

Finally, fix $b_0$.
The parameter spaces for positive multiples of $F_b-F_{b_0}$
to be rationally trivial have countably many constructible images in $B^{\circ}$.
None can be dense: a dominating parameter space would give~\eqref{eq:spread} with $Z=F_{b_0}$ and the same contradiction.
Their closures are therefore proper. Outside their countable union the difference is nonzero over $\mathbb{Q}$, as required.
\end{proof}

The fibre differences in Theorem~\ref{thm:fibres} have integral coefficients;
nonvanishing after tensoring with $\mathbb{Q}$ shows that they are not torsion.

\section{Canonical forms on universal Picard varieties}\label{sec:picard}

Let $M_{10}$ be the coarse moduli space of smooth curves of genus-$10$,
and let $S\subset M_{10}$ be its open subset parameterizing curves with trivial automorphism group.
There exists a universal curve
$\pi:\mathcal{C}\longrightarrow S.$
For $e\in\mathbb{Z}$, let
$$
p_e:P_e=\operatorname{Pic}^{e}(\mathcal{C}/S)\longrightarrow S
$$
be the relative Picard space.
Its geometric points are pairs $(C,L)$ with $\deg L=e$.
It is a smooth proper torsor over $S$ under the relative Jacobian $P_0$.
It is irreducible, and
$
\dim P_e=\dim S+10=27+10=37.
$
Write $B_e$ for  a smooth projective model of $P_e$.
Regular differential forms on smooth projective varieties are birational invariants,
so the spaces used below are independent of the choice of $B_e$.

\subsection{A canonical form in degree 10}

We first construct a form on $B_{10}$.
This requires a section of the canonical bundle itself, not just a section of a positive tensor power.

Let $\overline{\mathcal{M}}_{g,n}$ be the Deligne--Mumford stack of stable $n$-pointed genus-$g$ curves,
and let $\overline{M}_{g,n}$ be its coarse moduli space. Let
$$
\mathscr{M}=\overline{\mathcal{M}}_{10,10},\quad
\varphi:\mathscr{M}\longrightarrow\overline{\mathcal{M}}_{10},
$$
where $\varphi$ forgets the markings and stabilizes.
The stack $\mathscr{M}$ is smooth, proper, irreducible, and of dimension $37$.

We use additive notation in Picard groups:
for a Cartier divisor $D$, the symbol $[D]$ denotes the class of $\mathcal{O}(D)$.
The canonical line bundle of the smooth stack $\mathscr{M}$ is $\omega_{\mathscr{M}}=\det\Omega^1_{\mathscr{M}}$.
Let $\lambda_0$ denote the Hodge class on $\overline{\mathcal{M}}_{10}$ and set $\lambda=\varphi^*\lambda_0$.
Let $\psi_j$ be the first Chern class of the cotangent line at the $j$-th marking.
Write $\Delta_{\mathrm{irr}}$ for the reduced nonseparating boundary and $\delta_{\mathrm{irr}}$ for its class.
To index the separating boundary without repetition, let $I=\{1,\ldots,10\}$ and choose the pairs $(i,A)$ as follows:
$$
\begin{array}{ll}
i=0,&A\subseteq I,\quad |A|\geq2,\\
1\leq i\leq4,&A\subseteq I,\\
i=5,&A\subseteq I,\quad 1\in A.
\end{array}
$$
Their set is denoted by $\mathcal{I}$.
The divisor $\Delta_{i:A}$ has two components of genera $i$ and $10-i$ meeting at one node,
with the markings in $A$ on the genus-$i$ component; set $\delta_{i:A}=[\Delta_{i:A}]$.
The last convention chooses one representative of $(5,A)\sim(5,I\setminus A)$.

Let $\delta_0,\ldots,\delta_5$ be the reduced boundary classes on $\overline{\mathcal{M}}_{10}$, with $\delta_0$ nonseparating.
The standard pullback formulas are
\begin{equation}\label{eq:boundary-pullback}
\varphi^*(\delta_0)=\delta_{\mathrm{irr}},\quad
\varphi^*(\delta_i)=\sum_{(i,A)\in\mathcal{I}}\delta_{i:A}
\quad(1\leq i\leq5).
\end{equation}
All these are classes of reduced stack divisors, with multiplicity one in~\eqref{eq:boundary-pullback}.

Let $\mathscr{D}$ be the closure in $\mathscr{M}$ of the divisor
$
h^0\bigl(C,\mathcal{O}_C(x_1+\cdots+x_{10})\bigr)\geq2
$
on the smooth pointed locus.
Let $\mathscr{K}\subset\overline{\mathcal{M}}_{10}$ be the locus of curves lying on a $K3$ surface.
Set $a_{i,c}=\binom{|c-i|+1}{2}.$
We use the following known formulas for divisor classes, without repeating their proofs:
\begin{align}
[\mathscr{D}]&=-\lambda+\sum_{j=1}^{10}\psi_j
 -\sum_{(i,A)\in\mathcal{I}}a_{i,|A|}\delta_{i:A},
\label{eq:loganon}\\
[\mathscr{K}]&=7\lambda_0-\delta_0-\sum_{i=1}^{5}\beta_i\delta_i,
\quad (\beta_1,\beta_2,\beta_3,\beta_4)=(5,9,12,14),\quad\beta_5\geq6,
\label{eq:kthree}\\
[\omega_{\mathscr{M}}]&=13\lambda+\sum_{j=1}^{10}\psi_j
 -2\delta_{\mathrm{irr}}-2\sum_{(i,A)\in\mathcal{I}}\delta_{i:A}.
\label{eq:stackcanonical}
\end{align}
Formula~\eqref{eq:loganon} is Logan's calculation~\cite[Thm.~5.4]{Logan},
in the complete boundary form recalled in~\cite[P.~588]{FarkasVerra}.
Formula~\eqref{eq:kthree} is~\cite[Thm.~1.6]{FarkasPopa}.
Formula~\eqref{eq:stackcanonical} is the canonical class formula on the stack;
the coarse formula in~\cite[P.~592, Eq.~(1)]{FarkasVerra} has the additional term $-\delta_{1:\varnothing}$ from the elliptic tail involution.
We consistently use stack classes here.

The integral Picard groups of $\mathscr{M}$
and $\overline{\mathcal{M}}_{10}$ are free on the standard tautological classes,
subject only to the complementary boundary identifications.
This is the characteristic-zero theorem of Arbarello--Cornalba~\cite{ArbarelloCornalba}; see also~\cite[Thm.~B(1)]{FringuelliViviani}.
In particular, these groups have no torsion, and $\beta_5$ in~\eqref{eq:kthree} is an integer.
No value of $\beta_5$ more precise than the displayed lower bound is needed.

\begin{proposition}\label{prop:seed}
Every smooth projective model $B_{10}$ of $P_{10}$ satisfies
$$
H^0(B_{10},\Omega_{B_{10}}^{37})\neq0.
$$
\end{proposition}

\begin{proof}
For $(i,A)\in\mathcal{I}$ define
$$
c_{i,A}=\begin{cases}
\displaystyle\binom{|A|+1}{2}-2,&i=0,\\[4pt]
2\beta_i+a_{i,|A|}-2,&1\leq i\leq5.
\end{cases}
$$
These are nonnegative integers. Equations~\eqref{eq:boundary-pullback}--\eqref{eq:stackcanonical} give
\begin{equation}\label{eq:effective-canonical}
[\omega_{\mathscr{M}}]=
\left[\mathscr{D}+2\varphi^*\mathscr{K}
 +\sum_{(i,A)\in\mathcal{I}}c_{i,A}\Delta_{i:A}\right]
\end{equation}
in $\operatorname{Pic}(\mathscr{M})_{\mathbb{Q}}$.
Both sides are integral line bundle classes.
Torsion-freeness therefore upgrades~\eqref{eq:effective-canonical} to integral linear equivalence.
The effective divisor on the right is the zero divisor of a nonzero section $s\in H^0(\mathscr{M},\omega_{\mathscr{M}})$.

We check that $s$ is invariant under the symmetric group $\mathfrak{S}_{10}$ permuting the $10$ marked points.
Its divisor is invariant, so for each permutation $\sigma$,
the ratio $\sigma^*s/s$ is an invertible regular function on the connected proper stack $\mathscr{M}$, hence a constant.
These constants form a one-dimensional character.
Consider the action of the transposition $(12)$ at a general point of $\Delta_{0:\{1,2\}}$.
The rational tail has precisely two marked points and its attaching node.
Its unique automorphism exchanging the two markings acts as $-1$ on the node-smoothing direction and trivially on the boundary.
The general labelled curve at this point has trivial automorphism group.

Neither $\mathscr{D}$ nor $\varphi^*\mathscr{K}$ contains this boundary divisor as a component.
Since $c_{0,\{1,2\}}=1$,
local analytic coordinates can be chosen so that the transposition sends $(t,z_2,\ldots,z_{37})$ to $(-t,z_2,\ldots,z_{37})$ and
$$
s=t\,u(t,z)\,dt\wedge dz_2\wedge\cdots\wedge dz_{37},
\quad u(0,z)\neq0.
$$
The leading term is fixed: the two minus signs cancel.
Thus the character is trivial on a transposition, and hence on $\mathfrak{S}_{10}$.

Let $Q=\overline{M}_{10,10}/\mathfrak{S}_{10}$.
The invariant section descends to a rational top-degree form $\eta$ on $Q$.
It is regular on $Q_{\mathrm{reg}}$.
Indeed, on a local smooth finite quotient chart, a pole of order $a\geq1$
along a prime divisor downstairs would pull back with order $e-1-ea<0$,
where $e$ is the ramification index. This contradicts regularity of $s$.
Extension across codimension at least two on $Q_{\mathrm{reg}}$
then follows because its canonical sheaf is locally free.

By~\cite[Thm.~1.1]{FarkasVerra}, applied with $g=n=10$ and with the first tensor power,
$\eta$ extends to every resolution of $Q$. The Abel map
$$
Q\dashrightarrow P_{10},\quad
(C,x_1+\cdots+x_{10})\longmapsto
(C,\mathcal{O}_C(x_1+\cdots+x_{10}))
$$
is birational, as recalled in~\cite[P.~587]{FarkasVerra}:
on a general genus-ten curve,
a general degree-ten line bundle has a unique effective divisor.
Restriction to the automorphism-free locus does not change the function field.
Birational invariance of regular forms proves the proposition.
\end{proof}

The existence asserted here is stronger than the existence of a pluricanonical section.
In particular, it has not been inferred solely from the Kodaira dimension calculation of~\cite{FarkasVerra}.

\subsection{Changing the Picard degree}

Let
$$
\Lambda=\det(\pi_*(\omega_{\mathcal{C}/S})),\quad
\mathcal{A}=\omega_S\otimes\Lambda.
$$
Here $\omega_{\mathcal{C}/S}$ is the relative canonical line bundle.
Invariant differentials on the relative Jacobian and translation on its torsors give
\begin{equation}\label{eq:relative-canonical}
\omega_{P_e/S}\cong p_e^*(\Lambda),\quad
\omega_{P_e}\cong p_e^*(\mathcal{A}),\quad
(p_e)_*\mathcal{O}_{P_e}=\mathcal{O}_S.
\end{equation}
The last equality follows from properness and geometrically connected fibres.
Thus global top forms on $P_e$ correspond to sections of $\mathcal{A}$ on $S$.
These statements concern the relative Picard space and do not require a universal Poincar\'e line bundle on $\mathcal{C}\times_S P_e$.

\begin{proposition}\label{prop:all-degrees}
For every integer $e$, a smooth projective model $B_e$ of $P_e$ satisfies
$$
H^0(B_e,\Omega_{B_e}^{37})\neq0.
$$
\end{proposition}

\begin{proof}
Consider the morphism over $S$
\begin{equation*}
u:P_{10}\longrightarrow P_0,\quad
(C,L)\longmapsto(C,L^{\otimes9}\otimes\omega_C^{\otimes(-5)}).
\end{equation*}
It is well-defined since $9\cdot10-5\cdot18=0$.
On each geometric Picard fibre it is a translate of multiplication by $9$;
hence it is finite \mbox{\'etale}, of degree $9^{20}$.

Choose a smooth projective compactification of $P_{10}$ and a nonzero form $\Omega$ supplied by Proposition~\ref{prop:seed}. By~\eqref{eq:relative-canonical},
its restriction corresponds to a nonzero section $a\in H^0(S,\mathcal{A})$.
In the same identification for $P_0$, the pullback by $u$ multiplies a top relative differential by $9^{10}$,
while it is the identity on the base.
Therefore, there exists a regular top form $\eta$ on $P_0$, corresponding to $9^{-10}a$,
such that
\begin{equation}\label{eq:trace-nonzero}
u^*\eta=\Omega|_{P_{10}},\quad
\operatorname{Tr}_{u}(\Omega|_{P_{10}})=9^{20}\eta\neq0.
\end{equation}
Resolve the induced rational map between smooth projective models.
For a proper generically finite morphism of smooth varieties of equal dimension,
the trace map sends regular top forms to regular top forms; this is the standard duality trace   (see~\cite[Ex.~48.20.11 and Sec.~48.27]{Stacks}).
Its restriction to the open Picard space is the trace in~\eqref{eq:trace-nonzero}.
Consequently, $B_0$ has a nonzero regular $37$-form.
This also establishes extension across the boundary; regularity only on $P_0$ would not suffice.

For arbitrary $e$, the morphism
\begin{equation*}
v_e:P_e\longrightarrow P_0,\quad
(C,L)\longmapsto(C,L^{\otimes18}\otimes\omega_C^{\otimes(-e)})
\end{equation*}
is finite \mbox{\'etale} and surjective.
Indeed, its degree is zero on each curve,
and its fibrewise linear part is multiplication by $18$.
Resolve the rational map $B_e\dashrightarrow B_0$
and pull back the nonzero $37$-form.
Dominance gives nonvanishing; birational invariance descends the form from the resolution to $B_e$.
\end{proof}

\begin{remark}
If $3\nmid e$, integers $a,b$ can be chosen with $ae+18b=10$.
The map $(C,L)\mapsto(C,L^{\otimes a}\otimes\omega_C^{\otimes b})$
then gives a direct dominant morphism $P_e\longrightarrow P_{10}$.
This shorter argument does not cover $3\mid e$.
Proposition~\ref{prop:all-degrees} covers these remaining degrees because the trace in~\eqref{eq:trace-nonzero} is nonzero.
A trace of an arbitrary nonzero form can vanish; the equality $u^*\eta=\Omega|_{P_{10}}$,
obtained from the invariant relative differentials, is the necessary additional justification.
\end{remark}

The results of \cite{CasalainaMartinKassViviani} on the Kodaira dimension in all degrees
for compactified universal Jacobians provide related birational information.
We do not use them as a substitute for the canonical form argument above.

\section{A component of nonspecial space curves}\label{sec:components}

We next construct a component of $X_d$ that dominates $P_d$.
The irreducibility of the moduli space of smooth curves is due to Deligne--Mumford~\cite{DeligneMumford}.
Together with the geometrically irreducible Picard fibres, it gives the irreducibility used below.

A line bundle $L$ on a smooth curve $C$ is nonspecial if $H^1(C,L)=0$.
Write $\operatorname{Pic}^d(C)$ for the variety of degree-$d$ line bundles,
and $C^{(a)}$ for its $a$-th symmetric product, parameterizing effective divisors of degree $a$, with $C^{(0)}$ a point.
The following elementary dimension estimate records the precise degree bound needed in the construction.

\begin{lemma}\label{lem:line-bundles}
Let $C$ be a smooth projective curve of genus $10$ and let $d\geq13$.
Then a general $L\in\operatorname{Pic}^d(C)$ is nonspecial and very ample, and
$$
h^0(C,L)=d-9.
$$
A general four-dimensional subspace of $H^0(C,L)$ defines a closed embedding
$C\hookrightarrow\mathbb{P}^{3}$.
\end{lemma}

\begin{proof}
For $13\leq d\leq20$, consider the proper morphism
$$
C^{(2)}\times C^{(20-d)}\longrightarrow\operatorname{Pic}^d(C),
\quad (Z,A)\longmapsto\omega_C\otimes\mathcal{O}_C(Z-A).
$$
Its image has dimension at most $22-d\leq9$,
whereas $\dim\operatorname{Pic}^d(C)=10$.
If $L$ is outside this image, then
$$
H^0(C,\omega_C\otimes L^{-1}(Z))=0
$$
for every effective divisor $Z$ of degree two.
By Serre duality, $H^1(C,L(-Z))=0$ for every such $Z$, and also $H^1(C,L)=0$.
It follows that $H^0(C,L)\longrightarrow H^0(Z,L|_Z)$ is surjective for every length-two subscheme,  including $Z=2x$.
The standard very-ampleness criterion and Riemann--Roch~\cite[Chap.~IV]{Hartshorne} give the assertion.
For $d\geq21$ the same vanishing holds for every $L$, since $\deg(\omega_C\otimes L^{-1}(Z))=20-d<0$.

Set $r=h^0(C,L)=d-9$. If $r=4$, the complete linear series already embeds $C$ in $\mathbb{P}^{3}$.
If $r>4$, it embeds $C$ in $\mathbb{P}^{r-1}$. The closure of its secant and tangent lines has dimension at most three.
A general center $\mathbb{P}^{r-5}$ avoids this locus, since $(r-5)+3<r-1$, so projection is an embedding into $\mathbb{P}^{3}$.
These projections correspond to general four-dimensional subspaces of $H^0(C,L)$.
\end{proof}

\begin{proposition}\label{prop:component}
For every $d\geq13$, there exists an irreducible component $V_d\subset X_d$ of dimension $4d$ and a dominant rational map
\begin{equation}\label{eq:component-map}
\rho_d:V_d\dashrightarrow P_d,\quad
[C\hookrightarrow\mathbb{P}^{3}]\longmapsto(C,\mathcal{O}_C(1)).
\end{equation}
Its general geometric fibre is rational of dimension $q_d=4d-37$.
\end{proposition}

\begin{proof}
Let $P_d^{\circ}\subset P_d$ be the open locus of nonspecial very ample line bundles.
By Lemma~\ref{lem:line-bundles}, it is nonempty.
Fix homogeneous coordinates on $\mathbb{P}^{3}$.
Over a geometric point $(C,L)$,
the space of four ordered sections, modulo a common nonzero scalar, is
$\mathbb{P}\!\left(\operatorname{Hom}(\mathbb{C}^{4},H^0(C,L))\right).$
Here projectivization parameterizes one-dimensional subspaces of $\operatorname{Hom}(\mathbb{C}^{4},H^0(C,L))$.
The embedding condition defines a nonempty open subset of dimension
\begin{equation}\label{eq:fibre-dimension}
4h^0(C,L)-1=4(d-9)-1=q_d.
\end{equation}

The universal line bundle exists \mbox{\'etale} locally on $P_d^{\circ}$.
On overlaps, changing it by a line bundle from the base leaves the projectivized section space unchanged.
Thus the embedding loci descend to an algebraic space $T_d\longrightarrow P_d^{\circ}$
which is \mbox{\'etale} locally an open subset of a projective space bundle.
Globally, the ambient fibration may be a Brauer--Severi variety.
Its geometric fibres are irreducible, hence $T_d$ is irreducible and
\begin{equation*}
\dim T_d=37+q_d=4d.
\end{equation*}
The descended universal embedded curve gives a morphism to $\operatorname{Hilb}^{dt-9}(\mathbb{P}^{3})$,
the Hilbert scheme with Hilbert polynomial $dt-9$.
The image determines $C$, its hyperplane line bundle,
and its four coordinate sections up to common scalar.
Since curves in $S$ have trivial automorphism group, this morphism is injective on geometric points,
and its image has dimension $4d$.

We check the component assertion.
For a member $C\subset\mathbb{P}^{3}$,
let $L=\mathcal{O}_C(1)$.
The restricted Euler sequence and the normal sequence give
$$
\begin{gathered}
0\longrightarrow\mathcal{O}_C\longrightarrow L^{\oplus4}
\longrightarrow T_{\mathbb{P}^{3}}|_C\longrightarrow0,\\
0\longrightarrow T_C\longrightarrow T_{\mathbb{P}^{3}}|_C
\longrightarrow N_{C/\mathbb{P}^{3}}\longrightarrow0.
\end{gathered}
$$
Since $H^1(C,L)=0$, these imply
$$
H^1(C,T_{\mathbb{P}^{3}}|_C)=H^1(C,N_{C/\mathbb{P}^{3}})=0.
$$
Here $T_C$ and $T_{\mathbb{P}^{3}}$ are tangent bundles and $N_{C/\mathbb{P}^{3}}$ is the normal bundle.
The latter has rank two and degree $4d+18$, so, with $\chi=h^0-h^1$,
$$
\chi(C,N_{C/\mathbb{P}^{3}})=(4d+18)+2(1-10)=4d.
$$
Embedded deformation theory~\cite{Sernesi} shows that the Hilbert scheme is smooth of dimension $4d$ at $[C]$.
Consequently the irreducible closure of the image of $T_d$ is a Hilbert component.

Around a smooth reduced curve in characteristic zero, the Hilbert--Chow morphism is governed by the normal family comparison of Rydh~\cite[Paper~IV, Cor.~12.9]{Rydh}, as stated in~\cite[Rem.~2.1]{Ricolfi}.
The statement applies to a smooth projective ambient variety, without a Calabi--Yau assumption.
It follows that the image of our Hilbert component is an irreducible component $V_d$ of the reduced Chow variety $X_d$.

Recovering $(C,\mathcal{O}_C(1))$ on the dense smooth locus gives~\eqref{eq:component-map},
whose image contains $P_d^{\circ}$.
Its general geometric fibre is birational to the projective space in~\eqref{eq:fibre-dimension}.
This proves the proposition.
\end{proof}

\begin{remark}
The target $\mathbb{P}^{3}$ is fixed, so its projective coordinates are not quotiented out.
Equivalently, the fibre dimension is
$$
\dim\operatorname{Gr}(4,H^0(C,L))+\dim\operatorname{PGL}_4
=4(h^0(C,L)-4)+15=4h^0(C,L)-1.
$$
Here $\operatorname{Gr}(4,H^0(C,L))$ is the Grassmannian of four-dimensional subspaces
and $\operatorname{PGL}_4$ is the projective linear group. At $d=13$,
the geometric embedding fibres are $\operatorname{PGL}_4$.
Neither a trivialization of the corresponding torsor over the function field nor a splitting of the Brauer class is required.
\end{remark}

\section{Applications to Chow varieties}\label{sec:applications}

\subsection{Algebraically trivial cycles and odd homology}

\begin{proof}[Proof of Theorem~\ref{thm:chow}]
Let $d\geq13$.
Take the component $V_d$ of Proposition~\ref{prop:component}
and a smooth projective model $B_d$ of $P_d$. Proposition~\ref{prop:all-degrees} gives
$$
\dim V_d=4d,\quad \dim B_d=37,\quad
H^0(B_d,\Omega_{B_d}^{37})\neq0.
$$
Resolve the dominant rational map $V_d\dashrightarrow B_d$ to obtain
$$
f_d:Y_d\longrightarrow B_d,\quad
j_d:Y_d\longrightarrow V_d\hookrightarrow X_d.
$$
Theorem~\ref{thm:fibres}, applied with $N=4d$ and $m=37$,
proves the assertion of uncountable dimensionality .
The effective cycles are the proper pushforwards of the fundamental cycles of the fibres over a flat open subset $B_d^{\circ}\subset B_d$.
Their differences are
\begin{equation}\label{eq:alpha}
\alpha_{b,b_0}=(\Gamma_d)_*([b]-[b_0]),\quad
\Gamma_d=(f_d,j_d)_*[Y_d]\in\operatorname{CH}_{4d}(B_d\times X_d).
\end{equation}
They are integral classes, and Theorem~\ref{thm:fibres} shows that for every fixed $b_0\in B_d^{\circ}(\mathbb{C})$,
a very general difference is nonzero after tensoring with $\mathbb{Q}$. It is therefore not torsion.
\end{proof}

\begin{proof}[Proof of Theorem~\ref{thm:odd}]
Use the same component $V_d$ and a smooth projective resolution $Y_d$ of the graph of $V_d\dashrightarrow B_d$.
By Proposition~\ref{prop:all-degrees}, $B_d$ has a nonzero regular $37$-form.
 Dominant pullback yields $h^{37,0}(Y_d)\geq1$. Corollary~\ref{cor:general-hodge}, with $N=4d$ and $r=37$,
 gives the first inequality in~\eqref{eq:main-hodge}.
 Hodge symmetry gives the second.
 The two bidegrees are distinct, so $b_{k_d}(X_d)\geq2$.
 The equality with the integral homology rank follows from the universal coefficient theorem.
\end{proof}

\begin{proof}[Proof of Corollary~\ref{cor:24}]
Using the values $(q_d,k_d)=(15,67)$ for $d=13$ and $(59,155)$ for $d=24$, Theorems~\ref{thm:chow} and~\ref{thm:odd} give that the nonvanishing of the algebraically trivial subgroup implies that the two equivalence relations are distinct.
Moreover, Theorem~\ref{thm:odd} gives $b_{67}(X_{13})\geq2$ and $b_{155}(X_{24})\geq2$.
\end{proof}

\subsection{The degree-24 construction}

We describe the family underlying Corollary~\ref{cor:24}.
For $d=24$, a degree-$24$ line bundle on a smooth genus-$10$ curve is nonspecial and very ample,
and has $15$ sections.
For a general pair $(C,L)$, the embedding locus is a nonempty open subset of $\mathbb{P}^{59}$.
Let $F_{C,L}$ be the closure of its image in $X_{24}$.
The relevant dimensions are
$$
\dim B_{24}=37,\quad\dim V_{24}=96,\quad\dim F_{C,L}=59.
$$
On a common flat open subset of the Picard base, fix $(C_0,L_0)$ and choose $(C,L)$ very generally.
Then
$$
[F_{C,L}]-[F_{C_0,L_0}]
\in\operatorname{CH}_{59}(X_{24})_{\mathrm{alg}}
$$
is nonzero after tensoring with $\mathbb{Q}$.
These closure cycles agree with the push-forwards of the general fibre cycles of~\eqref{eq:alpha}:
a boundary component of a graph resolution
that dominates the base has smaller relative dimension,
and hence contributes no additional component of dimension 59 over a general point of the base.

The cycle class of this difference in $H_{118}(X_{24})$ is zero.
On the other hand, its family correspondence gives a nonzero action into $H_{155}(X_{24})$.
The cohomological Hodge types detected in weight $155$ are $(96,59)$ and $(59,96)$.

The whole Chow variety has dimension $327$, by the unrestricted dimension formula of Eisenbud--Harris~\cite[P.~295, Cor.]{EisenbudHarris}:
$$
\dim X_{24}=3+\frac{24(24+3)}2=327.
$$
This is not the dimension of $V_{24}$.
Since Proposition~\ref{prop:extreme} applies to irreducible components of any dimension,
the non-equidimensionality does not affect the argument.

The first few numerical instances are listed below.
The last column is the degree detecting the family correspondence, not the degree $2q_d$ of the cycle class of a fibre difference.
\begin{center}
\setlength{\tabcolsep}{12pt}
\begin{tabular}{@{}ccccc@{}}
\toprule
$d$ & $h^0(C,L)$ & $\dim V_d$ & $q_d$ & $k_d$\\
\midrule
13&4&52&15&67\\
14&5&56&19&75\\
18&9&72&35&107\\
24&15&96&59&155\\
\bottomrule
\end{tabular}
\end{center}
We make no minimality claim for $d$, $q_d$, or $k_d$.
The construction provides algebraic parameter spaces and fibre correspondences,
but it does not supply explicit polynomial equations for two prescribed representatives.
It establishes the existence of non-torsion classes.

\subsection{Comparison conjectures and degree stabilization}

\begin{proof}[Proof of Corollary~\ref{cor:comparisons}]
Theorem~\ref{thm:chow} gives an uncountably dimensional subspace of algebraically trivial classes. Such classes are homologically trivial. The singular homology of a complex projective variety is finite-dimensional over $\mathbb{Q}$, since the variety has the homotopy type of a finite CW complex. Both conclusions follow.
\end{proof}

\begin{proof}[Proof of Corollary~\ref{cor:stabilization}]
The proper pushforward commutes with the cycle class map and therefore preserves homological triviality. Since \(q_d\leq e\), the rational cycle class map on \(X_e\) is injective in dimension \(q_d\) by~\cite[Cor.~2.4]{ChenHuChow}. Hence every push-forward of an algebraically trivial class vanishes. Theorem~\ref{thm:chow} supplies an uncountably dimensional subspace of the kernel of this map.
\end{proof}

Since $s_{d,q_d}=s_{q_d-1,q_d}\circ\cdots\circ s_{d,d+1}$,
Corollary~\ref{cor:stabilization} implies that
at least one of the consecutive addition maps with source degree between d and $q_d-1$
fails to be injective on the $q_d$-dimensional Chow groups.
For $d=24$, such a transition occurs between degrees 24 and 59.
This does not contradict~\cite[Thm.~2.2]{ChenHuChow},
whose stabilization isomorphisms are asserted only when the cycle dimension is at most the source degree.

\subsection{Hodge structures, affine pavings, and Lawson homology}

\begin{proof}[Proof of Corollary~\ref{cor:paving}]
The two Hodge types in~\eqref{eq:main-hodge} are off the diagonal,
so $H^{k_d}(X_d,\mathbb{Q})$ is not Hodge--Tate.
A finite affine paving yields homology in even degrees and of Tate type, by Borel--Moore localization.
The nonzero odd homology in Theorem~\ref{thm:odd} therefore rules out such a paving.
\end{proof}

The virtual Hodge polynomial of a variety $X$ is
$$
E(X;u,v)=\sum_{j,a,b}(-1)^j
h^{a,b}\bigl(H_c^j(X,\mathbb{Q})\bigr)u^av^b,
$$
where $H_c^j$ denotes cohomology with compact support.
Unlike an individual cohomology group, the polynomial can exhibit cancellation between degrees.
Thus~\eqref{eq:main-hodge} alone does not disprove the virtual-Hodge assertion in~\cite[Conj.~4.1]{Hu2021}.
We do not compute $E(X_d;u,v)$ here.

\begin{proof}[Proof of Corollary~\ref{cor:lawson}]
By~\cite[Prop.~4.8]{Hu2021}, one has
$$
L_1H_j\bigl(C_{p,d}(\mathbb{P}^{n})\bigr)
 \cong H_j\bigl(C_{p,d}(\mathbb{P}^{n}),\mathbb{Z}\bigr)
\quad(j\geq2).
$$
Applying this comparison with $p=1$, $n=3$, and $j=k_d$,
tensoring with $\mathbb{Q}$, and invoking Theorem~\ref{thm:odd} gives the assertion.
\end{proof}

The classes in~\eqref{eq:alpha} vanish modulo algebraic equivalence.
Thus they do not provide nonzero classes in the Griffiths group,
the quotient of homologically trivial cycles by algebraically trivial cycles.
In particular, our Chow group counterexamples do not by themselves refute the independent Lawson homology assertions in~\cite{Hu2021,ChenHuLawson}.

\emph{Disclosure of AI assistance.}
The authors used ChatGPT 6 Pro for assistance with manuscript preparation
and revision. They take responsibility for the text and the
mathematical arguments.

\end{document}